\documentstyle{amsppt}

\def\tild{\widetilde}
\define\ord{\operatorname{ord}}

\redefine\phi{\varphi}
\topmatter
\title On the product property of the pluricomplex Green function \endtitle
\author
Armen Edigarian
\endauthor

\address
\noindent
Instytut Matematyki\newline
Uniwersytet Jagiello\'nski\newline
Reymonta 4\newline
30-059 Krak\'ow, Poland
\endaddress

\email
umedigar\@cyf-kr.edu.pl
\endemail

\abstract
We prove that the pluricomplex Green function has the 
pro\-duct property 
$g_{D_1\times D_2}=\max\{ g_{D_1},g_{D_2}\}$ 
for any domains $D_1\subset\Bbb C^n$ and $D_2\subset\Bbb C^m$.
\endabstract
\endtopmatter
\document
Let $E$ denote the unit disc in $\Bbb C$.
For any domain $G\subset\Bbb C^n$ define
$$
g_D(a,z):=\inf\Sb \phi\in\Cal O(E,D) \\
\phi (0)=z \\
a\in \phi(E)\endSb
\Big\{\prod_{\lambda\in\phi^{-1}(a)}
|\lambda|^{\ord_\lambda (\phi-a)}\Big\},\quad a,z\in D,
$$
where $\Cal O(E,D)$ denotes the set of all holomorphic
mappings $E\to D$ and $\ord_\lambda(\phi-a)$ denotes  multiplicity
of $\phi-a$ at $\lambda$.

The function $g_D$ is proposed by Poletsky (cf\. \cite{Pol}) and 
is called the {\it pluricomplex Green function for\/} $D$. 
We have that (see \cite{Jar-Pfl1}, Chapter IV)
$$
g_D(a,z):=\inf\Sb \phi\in\Cal O(\overline{E},D) \\
\phi (0)=z \\
a\in \phi(E)\endSb
\Big\{\prod_{\lambda\in\phi^{-1}(a)}
|\lambda|^{\ord_\lambda (\phi-a)}\Big\},\quad a,z\in D.
\tag a
$$
Note that in the formula (a) we take only $\lambda\in\phi^{-1}(a)$ such that
$\lambda\in E$;

\noindent
(b) For any domains $D_1$, $D_2$ and any holomorphic
mapping $f:D_1\to D_2$ we have the following contractible
property: $g_{D_2}\big(f(z),f(w)\big)\le
g_{D_1}(z,w)$, $z,w\in D_1$.

\bigskip

The main result of the paper is the following product property.
\proclaim{Theorem 1}
Let $D_1\subset\Bbb C^n$ and $D_2\subset\Bbb C^m$  be domains.
Then
$$
\multline
g_{D_1\times D_2}\big((z_1,w_1),(z_2,w_2)\big)=\max\{ g_{D_1}(z_1,z_2),
g_{D_2}(w_1,w_2)\},\\
(z_1,w_1),(z_2,w_2)\in D_1\times D_2.
\endmultline
$$
\endproclaim

\remark{Remark}
The product property for $D_1\times D_2$ for the pluricomplex
Green function  in the case when $D_1$ or $D_2$
is pseudoconvex was proved in \cite{Jar-Pfl2}.
Note that in \cite{Jar-Pfl2} the authors used the description
of the pluricomplex Green function given by M\. Klimek.
\endremark

\demo{Proof}
The inequality ''$\ge$'' follows from the property (b). 
So, we have to prove ''$\le$''.

Let $(a_1,b_1),(a_2,b_2)\in D_1\times D_2$. Suppose that
$N>0$ is such that

$$
\max\{ g_{D_1}(a_1,a_2),g_{D_2}(b_1,b_2)\}< N.\tag 1
$$
It is sufficient to prove that 
$$
g_{D_1\times D_2}\big((a_1,b_1),(a_2,b_2)\big)<N.
$$
There are holomorphic mappings $\phi_1:\bar E\to D_1$ and
$\phi_2:\bar E\to D_2$ such that $\phi_1(0)=a_2$, $\phi_2(0)=b_2$,
$$
\prod_{\lambda\in\phi_1^{-1}(a_1)}|\lambda|^{\ord_{\lambda}(\phi_1-a_1)}
<N\quad\text{ and }\quad
\prod_{\lambda\in\phi_2^{-1}(b_1)}
|\lambda|^{\ord_{\lambda}(\phi_1-b_1)}<N.
$$
Note that if $\phi:\bar E\to D_1$ is a holomorphic mapping,
$\phi(0)=a_2$,
and $\sigma_1,\dots,\sigma_l$ are zeros of $\phi-a_1$ counted
with multiplicities, then for a mapping
$$
\tild\phi(\lambda):=(\phi-a_2)\frac
{\prod_{j=1}^l(\lambda-\tild\sigma_j)}
{\prod_{j=1}^l(\lambda-\sigma_j)}\frac{\prod_{j=1}^l \sigma_j}
{\prod_{j=1}^l\tild \sigma_j}+a_2,\tag 2
$$
where $\tild \sigma_j$ is close enough to $\sigma_j$, 
$j=1,\dots,l$, we have 
$\tild\phi\in \Cal O(\bar E,D_1)$.
So, we may assume that $\ord_\lambda(\phi_1-a_1)=1$
for any $\lambda\in\phi_1^{-1}(a_1)\cap E$ and
$\ord_\lambda(\phi_2-b_1)=1$
for any $\lambda\in\phi_2^{-1}(b_1)\cap E$.

Let $\{\zeta_1,\dots,\zeta_\nu\}=\phi_1^{-1}(a_1)\cap E$ and
$\{\xi_1,\dots,\xi_\mu\}=\phi_2^{-1}(b_1)\cap E$.
We may assume that 
$|\zeta_1|<|\zeta_2|<\dots<|\zeta_\nu|$ and
$|\xi_1|<|\xi_2|<\dots<|\xi_\nu|$,
$\nu$, $\mu$ are minimal.
Then
$$
|\zeta_1\dots\zeta_\nu|\ge N|\zeta_\nu|^\nu\quad\text{ and }
\quad |\xi_1\dots\xi_\mu|\ge N|\xi_\mu|^\mu.
$$
For, if $|\zeta_1\dots\zeta_\nu|<N|\zeta_\nu|^\nu$ then
we may consider the mapping $\tild\phi_1(\lambda):=\phi_1
(\zeta_\nu \lambda)$, and it  contradicts the minimality of $\nu$.

Assume that $|\zeta_1\dots\zeta_\nu|<|\xi_1\dots\xi_\mu|$.
Then we consider the mapping 
$$
\tild\phi_1(\lambda):=\phi_1(t\lambda),\quad\lambda\in\bar E,
$$
where 
$t:=\Big(\frac{|\zeta_1\dots\zeta_\nu|}{|\xi_1\dots\xi_\mu|}\Big)^{\frac1\nu}$.
Then $\left|\frac 1t\zeta_j\right|<1$, $j=1,\dots,\nu$, and
$$
\left|\Big(\frac {\zeta_1}t\Big)\dots\Big(\frac {\zeta_\nu}t\Big)\right|=
|\xi_1\dots\xi_\mu|.
$$
Hence, we may assume that $|\zeta_1\dots\zeta_\nu|=|\xi_1\dots\xi_\mu|$.
Take $\tild\phi_1(\lambda)=\phi_1(e^{i\theta}\lambda)$, where
$\theta$ is such that $e^{-i\theta}\zeta_1\dots e^{-i\theta}\zeta_\nu=
\xi_1\dots\xi_\mu$. 
So, we may assume that
$\zeta_1\dots\zeta_\nu=\xi_1\dots\xi_\mu$.

Put
$$
B_1(\lambda):=\prod_{j=1}^\nu\frac{\zeta_j-\lambda}{1-\bar\zeta_j\lambda}
\quad\text{ and }\quad
B_2(\lambda):=\prod_{j=1}^\mu\frac{\xi_j-\lambda}{1-\bar\xi_j\lambda},
\quad\lambda\in\bar E,
$$
and $C:=B_1(0)=B_2(0)$. We have that $|C|<N$.

Put $A_1:=\{\lambda:B_1'(\lambda)=0\}$, $A_2:=\{\lambda:B_2'(\lambda)=0\}$, 
$A:=B_1(A_1)\cup B_2(A_2)$, $\tild A_1:=B_1^{-1}(A)$,
and $\tild A_2:=B_2^{-1}(A)$.

Put 
$$
\tild B_1(\lambda):=\frac{B_1(\lambda)-B_1(0)}
{1-\overline{B_1}(0)B_1(\lambda)}.
$$
It holds 
$$
B_1(\lambda)=\frac{B_1(0)+\tild B_1(\lambda)}
{1+\overline{B_1}(0)\tild B_1(\lambda)}.
$$
Note that $\tild B_1$ is a Blaschke product (because, it is a proper
function), $\tild B_1(\lambda)=0$ iff $B_1(\lambda)=B_1(0)$, and
$\tild B_1'(\lambda)=0$ iff $B_1'(\lambda)=0$.
Changing very little zeros of the $\tild B_1$, we may assume that
it has no multiple zeros and $\tild B_1(0)=0$. Then we have a new 
Blaschke product 
$$
\widehat B_1(\lambda):=\frac{B_1(0)+\tild B_1(\lambda)}
{1+\overline{B_1}(0)B_1(\lambda)},
$$ 
such that $\widehat B_1(0)=B_1(0)$ and
there is no $\lambda_0\in E$
such that $\widehat B_1(\lambda_0)=\widehat B_1(0)$ and 
$\widehat B_1'(\lambda_0)\not=0$.
The zeros of the function $\widehat B_1$ are close enough to
zeros of the $B_1$. So, we may assume that 
$C\not\in B_1(A_1)$ (replacing $B_1$ by $\widehat B_1$ and $\phi_1$
by $\widehat\phi_1$, where $\widehat\phi_1$ is constructed using (2)).
In a similar way, we may assume that $C\not\in B_2(A_2)$.
Hence, $C\not\in B_1(A_1)\cup B_2(A_2)$, and, therefore,
$0\not\in\tild A_1\cup \tild A_2$.

There exists a covering $\pi:E\to E\setminus A$ such 
that $\pi(0)=C$. Let us show that $\pi$ is a Blaschke product,
i\.e\.
$$
\lim_{r\to1}\int_0^{2\pi}\log|\pi(re^{i\theta})|d\theta=0.
$$
%(cf\. \cite{Col-Loh}, remarks after Theorem 2\.15, \cite{Sei}). 
Suppose that there exists $\lim_{r\to1}\pi(re^{i\theta})=t$
for some $e^{i\theta}\in\partial E$. 
We want to prove that $|t|=1$ or $t\in A$. 
Assume that $t\in E\setminus A$. 
Then there exists neighborhood $V$
of $t$ such that $\pi^{-1}(V)=\cup_j V_j$, $V_{j_1}\cap V_{j_2}=\varnothing$,
$j_1\not=j_2$.
There exists $r_0$ such that $\pi(re^{i\theta})\in V$,
$r\in[r_0,1)$.
There exists $j_0$ such that $r_0e^{i\theta}\in V_{j_0}$.
Then $re^{i\theta}\in V_{j_0}$ for any $r\in[r_0,1)$.
Hence, 
$$
e^{i\theta}=
\lim_{r\to1}(\pi|_{V_{j_0}})^{-1}\circ\pi(re^{i\theta})=
\pi^{-1}(t)\cap V_{j_0}.
$$
It is a contradiction. So, $|t|=1$ or $t\in A$. Note that
$A$ is a finite set, hence $\lim_{r\to1}\pi(re^{i\theta})\in\partial E$
for a\.a\. $\theta\in[0,2\pi)$. Therefore, $\pi$ is an inner function.
Hence, $\pi=BS$, where $B$ is a Blaschke product and
$$
S(\lambda)=\exp\Big(-\int_0^{2\pi}\frac{\lambda+e^{i\theta}}
{\lambda-e^{i\theta}}d\sigma(\theta)\Big),
$$
where $\sigma\ge0$ on $\partial E$ is a singular measure with respect to
Lebegue measure and $S^{\ast}(\lambda)=0$ for $\sigma$-almost all $\lambda\in
\partial E$. 
We know that if exists  $\lim_{r\to1}\pi(re^{i\theta})=t$, where $|t|<1$,
then $t\in A$. Hence, $t\not=0$. So $\sigma=0$ and $S\equiv1$,
therefore, $\pi$ is a Blaschke product.

Note that $B_1|_{E\setminus\tild A_1}:E\setminus \tild A_1\to E\setminus A$
and $B_2|_{E\setminus\tild A_2}:E\setminus \tild A_2\to E\setminus A$
are finite coverings. Hence, there exist liftings
$\psi_1:E\to E\setminus\tild A_1$, and
$\psi_2:E\to E\setminus\tild A_2$ such that
$\pi=B_1\circ\psi_1=B_2\circ\psi_2$,
$\psi_1(0)=\psi_2(0)=0$. 

We have that $\pi$ is a Blaschke product.
So, there exists $r<1$ such that
$$
\log|\pi(0)|-\frac 1{2\pi}\int_0^{2\pi}
\log|\pi(re^{i\theta})|d\theta<\log N.\tag 3
$$

Let $\lambda_1,\dots,\lambda_s$ be solutions of the equation
$\pi(r\lambda)=0$ counted with multiplicity.
By Jensen formula (cf\. \cite{Rud}, Theorem 15\.18) 
the left hand side of (3) is equal to 
$\sum_{j=1}^s\log|\lambda_j|$, and, therefore,
$\sum_{j=1}^s\log|\lambda_j|<\log N$.

Put 
$$
\gamma(\lambda):=\big(\phi_1\circ\psi_1(r\lambda),
\phi_2\circ\psi_2(r\lambda)\big).
$$
Note that $\gamma:E\to D_1\times D_2$,
$\gamma(0)=(a_1,b_1)$, and $\gamma(\lambda_j)=(a_2,b_2)$.
Hence,
$$
g_{D_1\times D_2}\big((a_1,b_1),(a_2,b_2)\big)\le
\prod_{j=1}^s|\lambda_j|<N.
$$
\qed
\enddemo

\head{Acknowledgement}\endhead
I would like to thank Professors M\. Jarnicki, P\. Pflug, and W\. Zwonek
for helpful discussions and remarks.

%%%%%%%%%%%%%%%%%%%%%%% Literatura %%%%%%%%%%%%%%%%%%%%%%%%%%%%%%

\enddocument